\documentclass[11pt,english,a4paper]{smfart}
\usepackage[english]{babel}
\usepackage{xcolor}
\usepackage{pbsi}
\usepackage[T1]{fontenc}
\usepackage{stmaryrd}
\usepackage{mathabx}
\usepackage[mathscr]{euscript}
\usepackage[pagebackref,colorlinks=true,linkcolor=blue,citecolor=blue,urlcolor=red, hypertexnames=true]{hyperref}

 \usepackage[abbrev,backrefs,nobysame]{amsrefs}
 \usepackage{caption}
 \usepackage{enumerate}

\usepackage{amsmath,amssymb,bm,amsfonts}
\usepackage{nicematrix}
\usepackage[all]{xy}
\entrymodifiers={+!!<0pt,\fontdimen22\textfont2>}

\usepackage{mathrsfs}

\usepackage{graphicx}

\usepackage{color}

\newtheorem{theorem}{Theorem}[section]
\newtheorem{lemma}[theorem]{Lemma}

\newtheorem{corollary}[theorem]{Corollary}
\newtheorem{proposition}[theorem]{Proposition}

\newtheorem*{remark*}{\it Remark}

\newcommand{\flba}[2]{
\xymatrix@C15pt{#1\ar@{|->}[r]&#2}}
\newcommand{\flcourte}[2]{
\xymatrix@C12pt{#1\ar[r]&#2}}

{\theoremstyle{definition}}

{\theoremstyle{definition}}

\theoremstyle{definition}

\newtheorem{fact}[theorem]{Fact}

\newtheorem{problem}[theorem]{Open Problem}

{\theoremstyle{definition}\newtheorem{remark}[theorem]{Remark}}

\def\D{\ensuremath{\mathbb D}}

\def\Z{\ensuremath{\mathbb Z}}

\def\C{\ensuremath{\mathbb C}}

\def\b{{{\mathcal B}}}

\def\bth{\begin{theorem}}
\def\blm{\begin{lemma}}
\def\bpr{\begin{proposition}}
\def\bpf{\begin{proof}}
\def\epf{\end{proof}}
\def\epr{\end{proposition}}
\def\elm{\end{lemma}}
\def\eth{\end{theorem}}
\def\bco{\begin{corollary}}
\def\eco{\end{corollary}}
\def\be{\begin{enumerate}}
\def\ee{\end{enumerate}}
\def\bea{\begin{enumerate}[\rm (a)]}
\def\beun{\begin{enumerate}[\rm (1)]}
\def\bei{\begin{enumerate}[\rm (i)]}

\newcommand{\vertiii}[1]{{\left\vert\kern-0.25ex\left\vert\kern-0.25ex\left\vert #1 
    \right\vert\kern-0.25ex\right\vert\kern-0.25ex\right\vert}}

\newcommand{\gd}{G_{\delta }}

\newcommand{\bbx}{{\mathcal{B}}_{1}(X)}
\newcommand{\ppx}{{\mathcal{P}}_{1}(X)}
\newcommand{\ppl}{{\mathcal{P}}_{1}(\ell_2)}

\newcommand{\wot}{\texttt{WOT}}
\newcommand{\sott}{\texttt{SOT}\mbox{$_{*}$}}

\newcommand{\sot}{\texttt{SOT}}
\newcommand{\sote}{\texttt{SOT}\mbox{$^{*}$}}
\newcommand{\bx}{{\mathcal B}(X)}
\newcommand{\px}{{\mathcal P}(X)}

\newcommand{\lan}{\langle}
\newcommand{\ran}{\rangle}

\numberwithin{equation}{section}

\author[Valentin Gillet]{Valentin Gillet}
\address[Valentin Gillet]{Univ. Lille, CNRS, UMR 8524 - Laboratoire Paul Painlevé, F-59000 Lille, France}
\email{valentin.gillet@univ-lille.fr}

\thanks{This work was supported in part by the project COMOP of the French National Research Agency (grant ANR-24-CE40-0892-01) and by the Labex CEMPI (ANR-11-LABX-0007-01).}

\begin{document}

\title[Similar operator topologies on the space of positive contractions]{Similar operator topologies on the space of positive contractions}


\begin{abstract} 
In this article, we study the similarity of the Polish operator topologies \wot\,, \sot\,, \sott\, and \sote\, on the set of the positive contractions on $\ell_p$ with $p > 1$. Using the notion of norming vector for a positive operator, we prove that these topologies are similar on $\ppl$, that is they have the same dense sets in $\ppl$. In particular, these topologies will share the same comeager sets in $\ppl$. We then apply these results to the study of typical properties of positive contractions on $\ell_p$-spaces in the Baire category sense. In particular, we prove that a typical positive contraction $T \in (\ppl, \sot)$ has no eigenvalue. This stands in strong contrast to a result of Eisner and Mátrai, stating that the point spectrum of a typical contraction $T \in (\mathcal{B}_1(\ell_2), \sot)$ contains the whole unit disk. As a consequence of our results, we obtain that a typical positive contraction $T \in (\ppl, \wot)$ (resp. $T \in (\ppl, \sott)$) has no eigenvalue.
\end{abstract}

\maketitle

\par\bigskip
\section{Introduction}\label{Section introduction}

If $X$ is a complex Banach space, we denote by $\bx$ the space of all bounded operators on $X$ and by $\bbx$ the set of all contractions on $X$, that is, bounded operators of norm at most $1$.

\smallskip

If $X$ has a Schauder basis $(e_n)_{n \geq 0}$, a vector $x = \displaystyle \sum_{n \geq 0} x_n e_n$ of $X$ is said to be positive (written $x \geq 0$) whenever $x_n \geq 0$ for every $n \geq 0$. A bounded operator $T$ on $X$ is said to be positive (written $T \geq 0$) whenever $T x \geq 0$ for every $x \geq 0$. If $X = \ell_p$ with $p \geq 1$ and if $x = \displaystyle \sum_{n \geq 0} x_n e_n$ is a vector of $X$, we denote by $\lvert x \rvert$ the positive vector $\lvert x \rvert := \displaystyle \sum_{n \geq 0 } \lvert x_n \rvert e_n$, which also belongs to $X$. We call $\lvert x \rvert$ the modulus of the vector $x$. We denote by $\px$ the set of all positive operators on $X$ and by $\ppx$ the set of all positive contractions on $X$. We refer to \cite{Gill} for more details on positive operators.

\smallskip

If $\tau$ is a topology on $\mathcal{Y} = \bbx$ or $\mathcal{Y} = \ppx$ such that $(\mathcal{Y}, \tau)$ is Polish, we say that a property $(P)$ of operators on $X$ is typical for $\tau$ if the set $\{ T \in \mathcal{Y} : T \; \textrm{satisfies} \; (P) \}$ is comeager in $\mathcal{Y}$, that is contains a dense $G_\delta$ subset of $\mathcal{Y}$. In this article, we will be working with the Strong Operator Topology, the Strong* Operator Topology, the Weak Operator Topology and the Dual Strong Operator Topology. The Strong Operator Topology is just the pointwise convergence topology on $\bx$, while the Strong* Operator Topology is the pointwise convergence topology for operators and their adjoints. The Weak Operator Topology is the weak pointwise convergence topology on $\bx$. Finally the Dual Strong Operator Topology is the pointwise convergence topology for the adjoints on $\bx$. We will respectively write \sot, \sote\,, \wot\, and \sott\, for these topologies. We can summarize the convergence of a sequence for these topologies in the following way:  
$$
\left\{
\begin{array}{ll}
         T_i \underset{i}{\overset{\sot}{\longrightarrow}}T \; \iff \; T_i x \underset{i}{\overset{\lVert . \rVert}{\longrightarrow}}Tx \quad  \textrm{for every} \; x \in X,\\
        T_i \underset{i}{\overset{\sote}{\longrightarrow}}T \; \iff \; T_i \underset{i}{\overset{\sot}{\longrightarrow}}T \quad \textrm{and} \quad T_i^* \underset{i}{\overset{\sot}{\longrightarrow}}T^*, \\
        T_i \underset{i}{\overset{\wot}{\longrightarrow}}T \; \iff \; \lan y^*, T_i x\ran \underset{i}{\overset{}{\longrightarrow}} \lan y^*, Tx \ran \quad \textrm{for every}\; y^* \in X^* \; \textrm{and every} \; x \in X, \\
        T_i \underset{i}{\overset{\sott}{\longrightarrow}}T \; \iff \; T_i^* \underset{i}{\overset{\sot}{\longrightarrow}}T^*.
    \end{array}
\right.
$$
For more details on these topologies, see \cite{GMM1}, \cite{GMM2}, \cite{Con} and \cite{Gill}. Also recall that when $X$ has a basis, the space $(\ppx, \sot)$ is Polish, that the space $(\ppx, \sote)$ is Polish when $X^*$ is separable and that the spaces $(\ppx, \wot)$ and $(\ppx, \sott)$ are Polish when $X$ is reflexive. This comes from the fact that $\ppx$ is closed in $\bbx$ for each one of these topologies.  

\smallskip

Typical properties of Hilbert spaces contractions for the \sot,\, the \sote\, and the \wot\, topologies have been studied by Eisner in \cite{E} and by Eisner and Mátrai in \cite{EM}. It was proved in \cite{E} that a typical contraction $T \in (\mathcal{B}_1(\ell_2), \wot)$ is unitary, and it was proved in \cite{EM} that a typical contraction $T \in (\mathcal{B}_1(\ell_2), \sot)$ is unitarily equivalent to the infinite-dimensional backward unilateral shift operator on $\ell_2(\Z_+ \times \Z_+)$. In particular, a typical contraction $T \in (\mathcal{B}_1(\ell_2), \wot)$ (resp. $T \in (\mathcal{B}_1(\ell_2), \sot)$) has a non-trivial invariant subspace. Typical properties of $\ell_p$-spaces contractions for the \sot\, and the \sote\, topologies have been studied in depth in \cite{GMM}, \cite{GMM1} and \cite{GMM2}. In particular, it was proved that a typical contraction $T \in(\mathcal{B}_1(\ell_1), \sot)$ (resp. $T \in (\mathcal{B}_1(\ell_2), \sote)$) has a non-trivial invariant subspace. The study of typical properties of positive contractions on $\ell_p$-spaces for the \sot\, and\, the \sote\, topologies was initiated in \cite{Gill}, and different phenomena were exhibited there. Indeed, it was proved in \cite[Proposition 3.7]{Gill} that when $p > 1$, a typical positive contraction $T \in (\mathcal{P}_1(\ell_p), \sot)$ (resp. $T \in (\mathcal{P}_1(\ell_p), \sote)$) is not a co-isometry, whereas a typical contraction  $T \in (\mathcal{B}_1(\ell_2), \sot)$ is a co-isometry. In particular, it cannot be asserted as in \cite[Theorem 5.2]{EM} that a typical positive contraction  $T \in (\mathcal{P}_1(\ell_2), \sot)$ is unitarily equivalent to the infinite-dimensional backward unilateral shift operator on $\ell_2(\Z_+ \times \Z_+)$. The point spectrum of an \sot-typical (resp. \wot-typical) positive contraction on $\ell_2$ is still unknown, whereas a typical positive contraction $T \in (\mathcal{P}_1(\ell_2), \sote)$ is such that $2T$ and $2T^*$ are hypercyclic by \cite[Proposition 3.2]{Gill}, so the point spectrum of a typical positive contraction $T \in (\ppl, \sote)$ is empty. Finally by \cite[Corollary 5.3]{EM}, a typical contraction $T \in (\mathcal{B}_1(\ell_2), \sot)$ has the open unit disk as point spectrum.

\smallskip

The notion of similar operator topologies was studied in \cite{GMM2} by Grivaux, Matheron and Menet. One advantage of similar topologies is that these topologies have the same comeager sets. It was proved in \cite{GMM2} that when $p > 2$, the topologies \sote\, and \sot\, are similar on $\mathcal{B}_1(\ell_p)$ and the topologies \wot\, and \sott\, are similar on $\mathcal{B}_1(\ell_p)$, and that when $1 < p < 2$, the topologies \wot\, and \sot\, are similar on $\mathcal{B}_1(\ell_p)$ and the topologies \sott\, and \sote\, are similar on $\mathcal{B}_1(\ell_p)$. It was pointed out in \cite[Theorem 2.6]{Gill} that this is still the case in $\mathcal{P}_1(\ell_p)$, since \cite[Lemma 5.17]{GMM1} remains true for positive contractions and since the proof of \cite[Theorem 3.4]{GMM2} is true for positive contractions as long as \cite[Lemma 5.17]{GMM1} holds for positive contractions. The proof of \cite[Theorem 3.4]{GMM2} uses appropriately the notion of norming vector for an operator.

\smallskip

The aim of this article is to link the comeager sets of the topologies \wot, \sot, \sott\, and\, \sote\, on $\mathcal{P}_1(\ell_2)$ by showing that that these topologies are similar on $\ppl$. This is strikingly different from what happens on $\mathcal{B}_1(\ell_2)$, where the topologies \wot, \sot, \sott\, and \sote\, are far from being similar. Similarity of these topologies on $\ppl$ will allow us to deduce some properties of a typical positive contraction on $\ell_2$ for the \sot,\, the \wot\, and the \sott\, topologies. In particular, we will be able to identify the point spectrum of a typical positive contraction $T$ on $\ell_2$ for the \wot\, and the \sot\, topologies. This answers an open problem from \cite[Question 5.2]{Gill}.

\subsection{Notations}\label{notation}
The following notation will be used throughout this article.
\par\smallskip

\smallskip

- The open unit disk of $\C$ will be denoted by $\D$.

\smallskip
- If $X$ is a Banach space, we denote by $S_X$ the unit sphere of $X$. 

\smallskip
- The closed linear span of a family of vectors $(x_i)_{i \in I} \subseteq X$ will be written as $[x_i : i \in I]$. 

\smallskip
- Let $(e_n)_{n \geq 0}$ be the canonical basis of $\ell_p$ for $p \geq 1$. We denote by $E_N$ the subspace $[e_0,...,e_N ]$ for every $N \geq 0$ and by $F_N$ the subspace $ [e_j : j > N]$ for every $N \geq 0$.

\smallskip
- For every $N \geq 0$, we denote by $P_N$ the canonical projection onto $E_N$ and by $Q_N$ the canonical projection onto $F_N$. The biorthogonal functionals sequence associated to the basis $(e_n)_{n \geq 0}$ will be denoted by $(e_n^*)_{n \geq 0}$.

\smallskip
- If $T \in \bx$, we respectively write $\sigma(T), \, \sigma_{ap}(T), \, \sigma_p(T)$ and $\sigma_{ess}(T)$ for the spectrum, the approximate spectrum, the point spectrum and the essential spectrum of $T$.

\smallskip
- If $T$ is an operator on $\ell_2$, we write $t_{i,j} := \lan e_i, T e_j \ran$ for every $i,j \geq 0$ and if $u$ is a vector of $\ell_2$, we write $u_k := \lan e_k, u \ran$ for every $k \geq 0$.

\smallskip
- If $u$ is a vector of $\ell_p$, we denote by $\textrm{Supp}(u) := \{ j \geq 0 : \lan e_j , u \ran \ne 0\}$ the support of $u$.

\smallskip

\subsection{Main results} \label{mainresults}

Given two topologies $\tau $ and $\tau'$ on a topological space $\mathbf{Y}$, we say that these topologies are similar if they have the same dense sets in $\mathbf{Y}$. It is well-known (see for instance \cite[Lemma 2.1]{GMM2}) that two similar topologies have the same comeager sets in $\mathbf{Y}$, but the converse is not always true as mentioned in \cite[Remark 2.2]{GMM2}. The article \cite{GMM2} studies the similarity of the Polish operator topologies \wot, \sot, \sott\, and \sote\, in $\mathcal{B}_1(\ell_p)$ for $p > 1$. A main tool in this investigation is the use of continuity points of the identity map from $(\mathcal{B}_1(\ell_p), \tau)$ into $(\mathcal{B}_1(\ell_p), \tau')$, where $\tau \subseteq \tau'$ are two distinct topologies among the topologies \wot, \sot, \sott\, and \sote. 

One of the main results of the article \cite{GMM2} is the following one, which only concerns the cases $1 < p < 2$ and $p > 2$. 

\smallskip

\bth[{\cite[Theorem 3.4]{GMM2}}] \label{intromainthgmm}
    Let $X = \ell_p$ with $p > 1$. If $1 < p < 2$, the topologies $\emph{\wot}$ and $\emph{\sot}$ are similar on $\bbx$ and the topologies $\emph{\sott}$ and $\emph{\sote}$ are similar on $\bbx$. If $p > 2$, the topologies $\emph{\sot}$ and $\emph{\sote}$ are similar on $\bbx$ and the topologies $\emph{\wot}$ and $\emph{\sott}$ are similar on $\bbx$.
\eth

\smallskip

Regarding the case $p = 2$, the topologies \wot, \sot, \sott\, and \sote\, are far from being similar on $\mathcal{B}_1(\ell_2)$. Indeed, it was proved in \cite[Corollary 5.3]{EM} that a typical contraction $T \in (\mathcal{B}_1(\ell_2), \sot)$ is such that $\sigma_p(T) = \D$, whereas a typical contraction $(\mathcal{B}_1(\ell_2), \sote)$ is such that $\sigma_p(T) = \emptyset$ by \cite[Proposition 7.1]{GMM1}. The latter property comes from the fact that a typical contraction $T \in (\mathcal{B}_1(\ell_2), \sote)$ is such that $2T^*$ is hypercyclic (\cite[Proposition 2.3]{GMM}), and it is very specific to the \sote\, topology. It was also proved in \cite[Theorem 3.1]{EM} that a typical contraction $T \in (\mathcal{B}_1(\ell_2), \wot)$ is unitary, whereas a typical contraction $T \in (\mathcal{B}_1(\ell_2), \sot)$ is not invertible by \cite[Corollary 3.5]{GMM1}.

The proof of Theorem \ref{intromainthgmm} uses the following result, proved in \cite[Corollary 7.4]{GMM2}, and stated in a weaker form here.

\smallskip

\bpr \label{introproptoposimlem}
Let $\mathbf{Y}$ be a set and let $\tau$ and $\tau'$ be two topologies on $\mathbf{Y}$. Suppose that the topologies $\tau$ and $\tau'$ are Polish on $\mathbf{Y}$ and such that $\tau \subseteq \tau'$. Then the following assertions are equivalent.
\begin{enumerate}[(a)]
    \item The topologies $ \tau$ and $\tau'$ are similar on $\mathbf{Y}$;
    \item The set of all points of continuity of the identity map from $(\mathbf{Y}, \tau)$ into $(\mathbf{Y}, \tau')$ is $\tau'$-dense in $\mathbf{Y}$.
\end{enumerate}
\epr

In this article, we study the similarity of the operator topologies \wot, \sot, \sott\, and \sote\, on $\mathcal{P}_1(\ell_p)$ for $p > 1$ and especially for $p = 2$. One of our tool will also be the points of continuity of the identity map on $\mathcal{P}_1(\ell_p)$ and, more precisely, Proposition \ref{introproptoposimlem}. 

The following result will also be used many times in this article and was proved in \cite[Proposition 5.15]{GMM1} for contractions. It is mentioned in \cite[Subsection 2d.]{Gill} that this Proposition remains true for positive contractions, since \cite[Lemma 5.17]{GMM1} remains true for positive contractions (see for instance \cite[Lemma 2.7]{Gill}). Before stating it, we will be using the following convention. 

Let $X = \ell_p$ with $p > 1$. If $T$ is an operator on $X$, we identify the operator $P_N T P_N$ on $X$ with the operator $P_N T_{|E_N}$ on $E_N$ and if $A$ is an operator on $E_N$, we identify $A$ with the operator $P_N A P_N$ on $X$. We can now state the corresponding result. 

\smallskip

\bpr \label{introprop5.15}
Let $X = \ell_p$ with $p > 1$. Let $N \geq 0$ and let $A \in \mathcal{P}_1(E_N)$. For any $\varepsilon > 0$, there exists an operator $B \in \mathcal{P}_1(E_{2N+1})$ satisfying
\begin{enumerate}[(i)]
    \item $\lVert B \rVert = 1$;
    \item $B$ has a norming vector $u \in E_{2N+1}$ for which $\textrm{Supp}(u) = \textrm{Supp}(Bu) = \{ 0,..., 2N+1 \}$;
    \item $\lVert B P_N - A \rVert < \varepsilon$.
\end{enumerate}
Moreover, if $\lVert A \rVert$ is sufficiently close to $1$, then one may require that in fact
\begin{enumerate}[(i')]
    \item $\lVert B - A \rVert < \varepsilon$.
\end{enumerate}
\epr

Our main results in this article are the following. The first one is the exact analogue of Theorem \ref{intromainthgmm} but in the positive setting.

\smallskip

\bth \label{intromainthp>1}
    Let $X = \ell_p$ with $p > 1$. If $1 < p < 2$, the topologies $\emph{\wot}$ and $\emph{\sot}$ are similar on $\ppx$ and the topologies $\emph{\sott}$ and $\emph{\sote}$ are similar on $\ppx$. If $p > 2$, the topologies $\emph{\sot}$ and $\emph{\sote}$ are similar on $\ppx$ and the topologies $\emph{\wot}$ and $\emph{\sott}$ are similar on $\ppx$.
\eth

\smallskip

The second one concerns the case $p = 2$ and is not true for contractions on $\ell_2$. This result also highlights the fact that the notion of typical property in $\ppl$ is very different from the notion of typical property in $\mathcal{B}_1(\ell_2)$ for each one of the topologies \wot, \sot, \sott\, and \sote. 

\smallskip

\bth \label{introthtoposimsurl2}
Any Polish operator topology lying between the topologies \emph{\wot}\, and \emph{\sote}\, is similar to \emph{\wot}\, and \emph{\sote}\, in $\ppl$. In particular the topologies \emph{\wot}, \emph{\sott}, \emph{\sot}\, and \emph{\sote}\, are similar on $\ppl$.
\eth

\smallskip

Using Theorem \ref{introthtoposimsurl2}, we give a negative answer to Question $5.2$ in \cite{Gill}.

\smallskip

\bth \label{introthvpspos}
Let $\tau$ be any Polish topology on $\ppl$ lying between \emph{\wot}\, and \emph{\sote}. A typical positive contraction $T \in (\ppl, \tau)$ is such that $T$ and $T^*$ have no eigenvalue.
\eth

\smallskip

The existence of a non-trivial invariant subspace for a typical positive contraction on $\ell_2$ for the \sot\, and the \sote\, topologies was investigated in \cite{Gill}. Thanks to this investigation and thanks to Theorem \ref{introthtoposimsurl2}, we also obtain the following result.

\smallskip

\bth \label{introinvsubtyppos}
Let $\tau$ be any Polish topology on $\ppl$ lying between \emph{\wot}\, and \emph{\sote}. A typical positive contraction $T \in (\ppl, \tau)$ has a non-trivial invariant subspace.
\eth

\smallskip

In Section \ref{Section3}, we propose a partial description of the points of continuity of the identity maps $\mathbf{i}_{\sot, \sote} : (\ppl, \sot) \to (\ppl, \sote), \; \mathbf{i}_{\wot, \sott} : (\ppl, \wot) \to (\ppl, \sott) $ and $\mathbf{i}_{\wot, \sot} : (\ppl, \wot) \to (\ppl, \sot)$. The description of these sets of points of continuity is interesting to study typical properties of positive contractions since these sets are dense-$\gd$ in $\ppl$ for the corresponding topologies (see \cite[Corollary 2.10]{GMM2}). Although the full description of these sets in the non-positive setting is quite simple (see \cite[Proposition 2.11]{GMM2}), we will see that such a description in the positive setting is more delicate.

Our main results in Section \ref{Section3} are the following. The first one gives a partial description of the sets of points of continuity of the identity maps for the topologies $(\sot, \sote)$ and $(\wot, \sott)$ on $\ppl$.

\smallskip

\bth \label{th4.6introsec4}
    Let $\mathcal{M}$ be the class of every positive contractions $T$ on $\ell_2$ such that 
    \begin{itemize}
        \item $\lVert T \rVert = 1$;
        \item there exists a family $(u_r)_{r \in J}$ of norming vectors for $T^*$ indexed by a set $J \subseteq \Z_+$ such that $u_r \geq 0$, $\lVert u_r \rVert = 1$ for every $r \in J$ and $\displaystyle \bigcup_{r \in J} \textrm{Supp}(u_r) = \Z_+$.
    \end{itemize}
    Then every operator in $\mathcal{M}$ is a point of continuity of both identity maps $\mathbf{i}_{\emph{\wot}, \emph{\sott}}$ and $\mathbf{i}_{\emph{\sot}, \emph{\sote}}$ on $\ppl$. 
\eth

\smallskip

The second result concerns the topologies $(\wot, \sot)$ on $\ppl$ and will be obtained from Theorem \ref{th4.6introsec4}.

\smallskip

\bth \label{th4.7introsec4}
Let $\mathcal{M}'$ be the class of every positive contractions $T$ on $\ell_2$ such that 
    \begin{itemize}
        \item $\lVert T \rVert = 1$;
        \item there exists a family $(u_r)_{r \in J}$ of norming vectors for $T$ indexed by a set $J \subseteq \Z_+$ such that $u_r \geq 0$, $\lVert u_r \rVert = 1$ for every $r \in J$ and $\displaystyle \bigcup_{r \in J} \textrm{Supp}(u_r) = \Z_+$.
    \end{itemize}
    Then every operator in $\mathcal{M}'$ is a point of continuity of the identity map $\mathbf{i}_{\emph{\wot}, \emph{\sot}}$ on $\ppl$. 
\eth

\smallskip

We end Section \ref{Section3} by proving that the classes $\mathcal{M}$ and $\mathcal{M}'$ introduced in Theorems \ref{th4.6introsec4} and \ref{th4.7introsec4} are large classes of points of continuity but we also show that the elements of $\mathcal{M}$ (resp. $\mathcal{M}'$) are not all the points of continuity of the maps $\mathbf{i}_{\wot,\sott}$ and $\mathbf{i}_{\sot, \sote}$ (resp. $\mathbf{i}_{\wot,\sot}$) on $\ppl$.

\section{Similar operator topologies and applications to typical properties of positive contractions} \label{sectiontool}

\subsection{Points of continuity of the identity map, norming vector for a positive operator} \label{subsect2a}

The notion of points of continuity of the identity map on $\mathcal{B}_1(\ell_p)$ was used in \cite{GMM2} to study the similarity of the topologies \wot, \sot, \sott\, and \sote\, in $\mathcal{B}_1(\ell_p)$. In this article, Proposition \ref{introproptoposimlem} will also play an important role to prove Theorems \ref{intromainthp>1} and \ref{introthtoposimsurl2}.

\smallskip

Let $X = \ell_p$ with $p > 1$. If $\tau$ and $\tau'$ are two Polish topologies on $\ppx$, we denote by $\mathcal{C}(\tau, \tau')$ the points of continuity of the identity map $\mathbf{i}_{\tau, \tau'} : (\mathcal{P}_1(X), \tau) \to (\mathcal{P}_1(X), \tau') $. If $\tau, \tau' \in \{ \wot, \sott, \sot, \sote\}$ are such that $\tau \subseteq \tau'$, then $\mathcal{C}(\tau, \tau')$ is $\tau$-$\gd$ dense in $\ppx$ (see for instance \cite[Corollary 2.10]{GMM2}). 

\smallskip

The notion of norming vector for an operator was also very useful in \cite{GMM2} and will also play an important role in this article. For example, this notion appears in Proposition \ref{introprop5.15}.

Recall that a vector $x \ne 0$ of $X = \ell_p$ is a norming vector for an operator $T \in \bx$ if $\lVert T x \rVert = \lVert T \rVert \lVert x \rVert $. If $x$ is a norming vector for $T$ and if $T$ is a positive operator, then $\lvert x \rvert$ is also a norming vector for $T$ with the same support as $x$. Indeed if $x \in S_X$ is a norming vector for $T$, then  $\lVert T \rVert = \lVert T x \rVert = \lVert \lvert T x \rvert \rVert \leq \lVert T \lvert x \rvert \rVert \leq \lVert T \rVert$, and so $\lVert T \lvert x \rvert \rVert = \lVert T \rVert$. Here we used the fact that the norm is increasing on $\ell_p$, as well as the fact that the norm of a vector is the norm of its modulus. The set of all norming vectors for $T \in \bx$ will be written $\mathcal{N}(T)$. If $\mathcal{N}(T) \ne \emptyset$, we say that the operator $T$ attains its norm. Also notice that there exist positive operators on $\ell_p$ which do not attain their norm. For example if $a = (a_n)_{n \geq 0}$ is a sequence of real numbers in $(0,1)$ such that $a_n \underset{n\to \infty}{\longrightarrow} 1$, then the positive diagonal operator $\Delta_a$ associated to the sequence $(a_n)_{n \geq 0}$ does not attain its norm.

\smallskip

The following proposition is a well-known fact and will be useful for our study. For the convenience of the reader, we recall here the proof of this fact.

\smallskip

\bpr \label{propbasevn}
Let $X = \ell_2$, let $T \in \bx$ and let $u \in X$ be a norming vector for $T$. Then $T^* T u = \lVert T \rVert^2 \, u$. In particular if $\mathcal{N}(T) \ne \emptyset$, then $u+ v \in \mathcal{N}(T)$ and $\lambda u \in \mathcal{N}(T)$ for every $u,v \in \mathcal{N}(T)$ and for every $\lambda \in \C$.
\epr

\smallskip

\bpf
We can suppose that $\lVert u \rVert = 1$. 

We clearly have $\lVert T^* T u \rVert \leq \lVert T \rVert^2$. Moreover, we have 
\begin{align*}
    \lVert T \rVert^2 &= \lVert T u \rVert^2 = \lan T u , T u \ran = \lan T^* T u , u \ran \leq \lVert T^* T u \rVert 
\end{align*}
and thus $\lVert T^* T u \rVert = \lVert T \rVert^2$. By the equality case of the Cauchy–Schwarz inequality, there is a scalar $\alpha$ such that $T^* T u = \alpha u$, and we must have $\alpha = \lVert T \rVert ^2$. This proves Proposition \ref{propbasevn}.
\epf

\begin{remark} \label{remarkvectnormant}
    Conversely, if there exists a vector $u \in \ell_2$ such that $\lVert u \rVert = 1$ and $T^* T u = \lVert T \rVert^2 \, u$, then $u \in \mathcal{N}(T)$ since $\lVert T u \rVert ^2 = \lan T^* T u , u \ran = \lVert T \rVert^2 $. In particular if a vector $u \in \ell_2$ is norming for $T$, then $Tu$ is a norming vector for $T^*$.
\end{remark}

\subsection{Proof of Theorem \ref{intromainthp>1}}

The proof of Theorem \ref{intromainthgmm} uses the following proposition, which is a consequence of \cite[Proposition 5.15]{GMM1}. It is the analog of Proposition \ref{introprop5.15} but for contractions instead of positive contractions.  

\smallskip

\bpr \label{lem3.2gmm2}
Let $p > 1$ and let $X = \ell_p$. Let $\mathcal{U}$ be a non-empty \emph{\sote}-open set in $\mathcal{B}_1(X)$ and let $n_0 \geq 0$. Then one can find an index $M \geq n_0$ and an operator $B \in \mathcal{B}_1(E_M) $ such that $B \in \mathcal{U}$ and such that $B$ satisfies the following properties: $\lVert B \rVert = 1$ and $B$ has a norming vector $u \in E_M$ for which $\textrm{Supp}(Bu) = \{ 0,...,M \}$. 
\epr

\smallskip

The proof of Theorem \ref{intromainthgmm} also uses the following result, which follows from \cite[Lemma 3.3]{GMM2} and also works for positive contractions.

\smallskip

\bpr \label{lem3.3gmm2}
Assume that $p > 2$. Let $M \geq 0$ and let $B \in \mathcal{B}_1(E_M)$. Suppose that $\lVert B \rVert = 1$ and that $B$ has a norming vector $u \in E_M$ for which $\textrm{Supp}(Bu) = \{ 0,...,M \}$. Then for every $\varepsilon > 0$, there exists $\delta > 0$ such that the following property holds: if $T \in \mathcal{B}_1(\ell_p)$ is such that $\lVert P_M (T-B) P_M \rVert < \delta$, then $\lVert P_M (T-B) \rVert < \varepsilon$.
\epr

\smallskip

In order to adapt the proof of \cite[Theorem 3.4]{GMM2} in the positive setting, that is in order to prove Theorem \ref{intromainthp>1}, we only have to show that Proposition \ref{lem3.2gmm2} is working for positive contractions. This is the aim of the following proposition.

\smallskip

\bpr \label{lem3.2gmmpos}
Let $p > 1$ and let $X = \ell_p$. Let $\mathcal{U}$ be a non-empty \emph{\sote}-open set in $\mathcal{P}_1(X)$ and let $n_0 \geq 0$. Then one can find an index $M \geq n_0$ and a positive operator $B \in \mathcal{P}_1(E_M) $ such that $B \in \mathcal{U}$ and such that $B$ satisfies the following properties: $\lVert B \rVert = 1$ and $B$ has a norming vector $u \in E_M$ for which $\textrm{Supp}(Bu) = \{ 0,...,M \}$.
\epr

\smallskip

\bpf
Since the set $S^+(X) := \{ T \in \ppx : \lVert T \rVert = 1 \}$ is \sote-dense in $\ppx$ by \cite[Proposition 3.1]{Gill}, there exists a positive operator $A$ in\, $\mathcal{U} \cap S^+(X)$. Let $\varepsilon > 0$. There exists a vector $x \in S_X$ such that $\lVert A x \rVert > 1- \varepsilon$. Moreover, one can find an integer $n > n_0$ such that $A_n := P_n A P_n$ belongs to $\mathcal{U}$ and $\lVert A_n x \rVert > \lVert A x \rVert - \varepsilon > 1 - 2 \varepsilon$. Thus the norm of $A_n$ is sufficiently close to $1$ and, applying Proposition \ref{introprop5.15}, one can find an index $M > n_0$ and a positive operator $B \in \mathcal{P}_1(E_M)$ in $\mathcal{U}$ satisfying all the properties listed in Proposition \ref{lem3.2gmmpos}. This concludes the proof of Proposition \ref{lem3.2gmmpos}.
\epf

\smallskip

By Propositions \ref{lem3.3gmm2} and \ref{lem3.2gmmpos}, the proof of Theorem \ref{intromainthgmm} given in \cite[Theorem 3.4]{GMM2} also applies to positive contractions, and thus one can state the following theorem.

\smallskip

\bth \label{mainthmposp>1}
Let $X = \ell_p$ with $p > 1$. If $1 < p < 2$, the topologies $\emph{\wot}$ and $\emph{\sot}$ are similar on $\ppx$ and the topologies $\emph{\sott}$ and $\emph{\sote}$ are similar on $\ppx$. If $p > 2$, the topologies $\emph{\sot}$ and $\emph{\sote}$ are similar on $\ppx$ and the topologies $\emph{\wot}$ and $\emph{\sott}$ are similar on $\ppx$.
\eth

\subsection{Proof of Theorem \ref{introthtoposimsurl2}}

We now focus on the case $p = 2$. As mentioned in the introduction, the topologies \wot, \sot, \sott\, and \sote\, are far from being similar on $\mathcal{B}_1(\ell_2)$. It appears that the situation is considerably different for positive contractions as we will see. This is not surprising since a typical positive contraction $T \in (\ppl, \sot)$ is not a co-isometry by \cite[Proposition 3.7]{Gill}, and so we cannot apply the Wold decomposition to a typical positive contraction on $\ell_2$ as in \cite[Theorem 5.2]{EM}.

In order to prove Theorem \ref{introthtoposimsurl2}, we need to adapt Proposition \ref{lem3.3gmm2} to positive contractions on $\ell_2$. This is the aim of the next result.

\smallskip

\bpr \label{lem3.3pos2}
Let $M \geq 0$ and let $B \in \mathcal{P}_1(E_M)$. Suppose that $\lVert B \rVert = 1$ and that $B^*$ has a positive norming vector $u \in E_M$ for which $\textrm{Supp}(u) = \{ 0,...,M \}$. Then for every $\varepsilon > 0$, there exists $\delta > 0$ such that the following property holds: if $T \in \ppl$ is such that $\displaystyle \max_{0 \leq k,l \leq M} \lvert \lan e_k, (T-B) e_l \ran \rvert < \delta$, then $\displaystyle \max_{0 \leq k \leq M} \lVert (T-B)^* e_k \rVert < \varepsilon$.
\epr

\smallskip

\bpf
We can suppose that $\lVert u \rVert = 1$. We also set $u_l := \lan e_l, u \ran$ for every $0 \leq l \leq M$ and $b_{i,j} := \lan e_i, B e_j \ran$ for every $0 \leq i,j \leq M$. 
\smallskip

Let $\varepsilon > 0$ and let $T \in \ppl$ be such that $\displaystyle \max_{0 \leq k,l \leq M} \lvert \lan e_k, (T-B) e_l \ran \rvert < \delta$, where $\delta > 0$ is such that
\begin{align}
   & b_{i,j} - \delta > 0 \quad \textrm{for every} \quad 0 \leq i,j \leq M \quad \textrm{such that} \; b_{i,j} > 0 , \label{lem3.3eq1} \\
    & \textrm{and} \quad 2 \delta \displaystyle \max_{0 \leq l \leq M} \frac{1}{u_l^2} \left(\displaystyle \sum_{j=0}^M u_j^2 \sum_{k=0}^M b_{j,k} + \displaystyle \sum_{0 \leq i < j \leq M} u_i u_j \sum_{k=0}^M (b_{i,k} + b_{j,k}) \right) < \varepsilon. \label{lem3.3eq2}
\end{align}

Condition (\ref{lem3.3eq2}) makes sense since $u_l > 0$ for every $0 \leq l \leq M$, by the assumption on the support of the norming vector $u$.

\medskip

In order to prove that $\displaystyle \max_{0 \leq k \leq M} \lVert (T-B)^* e_k \rVert^2 < \varepsilon$, we first need to control each quantity $\lVert T^* e_l \rVert^2$, $0 \leq l \leq M$.

\smallskip

Using the parallelogram identity, one has
\begin{align}
    1 \geq \Big\lVert T^* \Big(\displaystyle \sum_{j=0}^M u_j e_j \Big) \Big\rVert^2 =  \displaystyle \sum_{j=0}^M u_j^2 \lVert T^* e_j \rVert^2 + 2 \displaystyle \sum_{0 \leq i < j \leq M} u_i u_j \, \lan T^* e_i, T^* e_j \ran. \label{lem3.3eq3}
\end{align}

Since the operator $T$ is positive, each quantity $\lan T^* e_i, T^* e_j \ran$, $0\leq i,j \leq M$, is non-negative. By positivity and since $\displaystyle \max_{0 \leq k,l \leq M} \lvert \lan e_k, (T-B) e_l \ran \rvert < \delta$, one also has 
\begin{align}
    \lan T^* e_i, T^* e_j \ran \geq \displaystyle \sum_{\substack{0 \leq k \leq M \\ b_{i,k} > 0, \, b_{j,k} > 0}} (b_{i,k} - \delta)( b_{j,k} - \delta) \label{lem3.3eq4}
\end{align}
for every $0 \leq i,j \leq M$.

\smallskip

Thus using (\ref{lem3.3eq3}), (\ref{lem3.3eq4}) and using the positivity of the vector $u$, we can control $\lVert T^* e_l \rVert^2$ for every $0 \leq l \leq M$ in the following way:
\begin{align}
    u_l^2 \lVert T^* e_l \rVert^2 \leq 1 - \displaystyle \sum_{\substack{0 \leq j \leq M \\ j \ne l }} u_j^2 \lVert T^* e_j \rVert^2 - 2 \displaystyle \sum_{0 \leq i<j \leq M} u_i u_j \sum_{\substack{0 \leq k \leq M \\ b_{i,k} > 0, \, b_{j,k} > 0}} (b_{i,k} - \delta)( b_{j,k} - \delta). \label{lem3.3eq5}
\end{align}

\smallskip

Exactly as (\ref{lem3.3eq4}), one has
\begin{align}
    \lVert T^*e_j \rVert^2 \geq \displaystyle \sum_{\substack{0 \leq k \leq M \\ b_{j,k} > 0}} ( b_{j,k} - \delta)^2\label{lem3.3eq6}
\end{align}
for every $0 \leq j \leq M$.

\smallskip

Thus using (\ref{lem3.3eq5}) and (\ref{lem3.3eq6}), one obtains for every $0 \leq l \leq M$
\begin{align}
   \label{lem3.3eq7} u_l^2 \lVert T^* e_l \rVert^2 &\leq 1 - \displaystyle \sum_{\substack{0 \leq j \leq M \\ j \ne l }} u_j^2 \displaystyle \sum_{\substack{0 \leq k \leq M \\ b_{j,k} > 0}} ( b_{j,k} - \delta)^2 \\ \notag
    &- 2 \displaystyle \sum_{0 \leq i<j \leq M} u_i u_j \sum_{\substack{0 \leq k \leq M \\ b_{i,k} > 0, \, b_{j,k} > 0}} (b_{i,k} - \delta)( b_{j,k} - \delta). 
\end{align}

\smallskip

In order to prove that $\displaystyle \max_{0 \leq k \leq M} \lVert (T-B)^* e_k \rVert^2 < \varepsilon$, we also need to control $\lan T^* e_l, B^* e_l \ran$ for every $0 \leq l \leq M$. But as (\ref{lem3.3eq4}),
\begin{align}
    \lan T^* e_l , B^* e_l \ran \geq \displaystyle \sum_{\substack{0 \leq k \leq M \\ b_{l,k} > 0}} b_{l,k} (b_{l,k} - \delta) \label{lem3.3eq8}
\end{align}
for every $0 \leq l \leq M$.

\smallskip
Thus combining (\ref{lem3.3eq7}) and (\ref{lem3.3eq8}), the expression 
\begin{align}
  \label{lem3.3eq9}  u_l^2 \lVert T^* e_l - B^* e_l \rVert^2 = u_l^2 \lVert T^* e_l \rVert^2 - 2 u_l^2 \lan T^* e_l , B^* e_l \ran + u_l^2 \lVert B^* e_l \rVert^2
\end{align}
is bounded above for every $0 \leq l \leq M$ by
\begin{align}
   \label{lem3.3eq10} 1 &- \displaystyle \sum_{\substack{0 \leq j \leq M \\ j \ne l }} u_j^2 \displaystyle \sum_{\substack{0 \leq k \leq M \\ b_{j,k} > 0}} ( b_{j,k} - \delta)^2 - 2 \displaystyle \sum_{0 \leq i<j \leq M} u_i u_j \sum_{\substack{0 \leq k \leq M \\ b_{i,k} > 0, \, b_{j,k} > 0}} (b_{i,k} - \delta)( b_{j,k} - \delta) \\ \notag
    &-2 u_l^2 \displaystyle \sum_{\substack{0 \leq k \leq M \\ b_{l,k} > 0}} b_{l,k} (b_{l,k} - \delta) + u_l^2 \displaystyle \sum_{0 \leq k \leq M} b_{l,k}^2.
\end{align}

\smallskip

Now putting together the terms in (\ref{lem3.3eq10}), the expression (\ref{lem3.3eq9}) is bounded above for every $0 \leq l \leq M$ by

\begin{align}
   \label{lem3.3eq11} & 1- \left(\displaystyle \sum_{j=0}^M u_j^2 \displaystyle \sum_{k=0}^M b_{j,k}^2 + 2 \displaystyle \sum_{0 \leq i < j \leq M} u_i u_j  \displaystyle \sum_{k=0}^M b_{i,k} b_{j,k}  \right)\\ \notag 
   &+ 2 \delta \left(\displaystyle \sum_{j=0}^M u_j^2 \sum_{k=0}^M b_{j,k} + \displaystyle \sum_{0  \leq i < j \leq M} u_i u_j \sum_{k=0}^M (b_{i,k} + b_{j,k}) \right).
\end{align}
Let us remark that
\begin{align*}
    \displaystyle \sum_{k=0}^M \lan u, B e_k \ran^2 &= \displaystyle \sum_{k=0}^M \left( \displaystyle \sum_{j=0}^M u_j b_{j,k} \right)^2 \\
    &=\displaystyle \sum_{j=0}^M u_j^2 \displaystyle \sum_{k=0}^M b_{j,k}^2 + 2 \displaystyle \sum_{0 \leq i < j \leq M} u_i u_j  \displaystyle \sum_{k=0}^M b_{i,k} b_{j,k}.
\end{align*}
Thus using (\ref{lem3.3eq11}), the expression (\ref{lem3.3eq9}) can be bounded above for every $0 \leq l \leq M$ by
\begin{align}
   \label{lem3.3eq12}  \left(1- \displaystyle \sum_{j=0}^M \lan u , B e_j \ran^2 \right) + 2 \delta \left(\displaystyle \sum_{j=0}^M u_j^2 \sum_{k=0}^M b_{j,k} + \displaystyle \sum_{0  \leq i < j \leq M} u_i u_j \sum_{k=0}^M (b_{i,k} + b_{j,k}) \right).
\end{align}

\smallskip
Now noticing that $$\displaystyle \sum_{j=0}^M \lan u , B e_j \ran^2 = \lVert B^* u \rVert^2 = 1, $$ 
it follows that the expression (\ref{lem3.3eq9}) is bounded above by
\begin{align}
    2 \delta \left(\displaystyle \sum_{j=0}^M u_j^2 \sum_{k=0}^M b_{j,k} + \displaystyle \sum_{0  \leq i < j \leq M} u_i u_j \sum_{k=0}^M (b_{i,k} + b_{j,k}) \right)
\end{align}
for every $0 \leq l \leq M$.

\smallskip

From (\ref{lem3.3eq1}) and (\ref{lem3.3eq2}), it easily follows that $\displaystyle \max_{0 \leq k \leq M} \lVert (T-B)^* e_k \rVert^2 < \varepsilon$. This concludes the proof of Proposition \ref{lem3.3pos2}.
\epf

Now with Propositions \ref{lem3.2gmmpos} and \ref{lem3.3pos2}, we are able to prove the following theorem regarding the case $p = 2$ in the positive setting.

\smallskip

\bth \label{bigthpossiml2}
Any Polish operator topology on $\ppl$ lying between the topologies \emph{\wot}\, and \emph{\sote}\, is similar to \emph{\wot}\, and \emph{\sote}\, in $\ppl$.
\eth

\smallskip

\bpf
In order to prove Theorem \ref{bigthpossiml2}, it is enough to prove that the topologies \wot\, and \sote\, are similar on $\ppl$. 

If we prove that the set $\mathcal{C}(\wot, \sott)$ is \sote-dense in $\ppl$, then the set $\mathcal{C}(\sot, \sote)$ will also be \sote-dense in $\ppl$, since $\mathcal{C}(\wot, \sott) \subseteq \mathcal{C}(\sot, \sote)$. Since the map $T \mapsto T^*$ is a homeomorphism from $(\mathcal{C}(\wot, \sott), \sote)$ onto $(\mathcal{C}(\wot, \sot), \sote)$ (see Proposition \ref{proprmkptscont} for more details), the set $\mathcal{C}(\wot, \sot)$ will also be \sote-dense in $\ppl$. Using Proposition \ref{introproptoposimlem}, and using the fact that the relation "being similar to" is a transitive relation on the set of Polish topologies on $\ppl$, the topologies \wot\, and \sote\, will be similar on $\ppl$.

Hence to prove Theorem \ref{bigthpossiml2}, it suffices to prove that $\mathcal{C}(\wot, \sott)$ is \sote-dense in $\ppl$.

It is not difficult to show that the metric
$$
d : (S, T) \mapsto \displaystyle \sum_{n=0}^{\infty} 2^{-n} \lVert (T-S)^* e_n \rVert
$$
generates the \sott\, topology on $\ppl$. Thus one has
\begin{align*}
\mathcal{C}(\wot, \sott) &= \displaystyle \bigcap_{\eta > 0} \mathcal{V}_{\eta} = \bigcap_{k \geq 1} \mathcal{V}_{1/k}\, ,
\end{align*}
where $$\mathcal{V}_{\eta} := \{ T \in \ppl : \textrm{there exists a \wot-neighborhood} \; \mathcal{W} \; \textrm{of} \; T \, \textrm{such that}\; diam(\mathcal{W}) < \eta \}$$ and where $diam(\mathcal{W})$ is the diameter of $\mathcal{W}$ for the distance $d$. Every set $\mathcal{V}_{\eta}$ is \wot-open and hence \sote-open in $\ppl$. Thus in order to prove that $\mathcal{C}(\wot, \sott)$ is \sote-dense in $\ppl$, it is enough to prove that each $\mathcal{V}_{\eta}$ is \sote-dense in $\ppl$ by the Baire category theorem.

Let $\mathcal{U}$ be a non-empty \sote-open set of $\ppl$ and let $n_0 \geq 0$ be such that $$\displaystyle \sum_{n > n_0} 2^{-n} < \eta/8.$$ By Proposition \ref{lem3.2gmmpos}, there exist $M > n_0$ and $B \in \mathcal{P}_1(E_M)$ such that $B \in \mathcal{U}$, $\lVert B \rVert = 1$ and $B$ has a positive norming vector $z \in E_M$ for which $\textrm{Supp}(Bz) = \{0,...,M \} $. Thus by Remark \ref{remarkvectnormant}, $B^*$ has a positive norming vector $u \in E_M$ whose support is $\{ 0,..., M \}$. By Proposition \ref{lem3.3pos2}, there exists $\delta > 0$ such that the following property holds: if $T \in \ppl$ is such that $\displaystyle \max_{0 \leq k,l \leq M} \lvert \lan e_k, (T-B) e_l \ran \rvert < \delta$, then $\displaystyle \max_{0 \leq k \leq M} \lVert (T-B)^* e_k \rVert < \eta / 8$. Now let us set $\mathcal{W} := \{ T \in \ppl :\displaystyle \max_{0 \leq k,l \leq M} \lvert \lan e_k, (T-B) e_l \ran \rvert < \delta\}$. The set $\mathcal{W}$ is a \wot-neighborhood of $B$ and $diam(\mathcal{W}) < \eta$ since for every $S,T \in \mathcal{W}$, we have
\begin{align*}
d(S,T) &\leq d(S,B) + d(B,T) \\
&\leq 2 \displaystyle \sum_{n=0}^M 2^{-n} \eta/8 + 4 \displaystyle \sum_{n > M}  2^{-n}  \\
&< \eta.
\end{align*}
This proves that $\mathcal{U} \cap \mathcal{V}_{\eta} \ne \emptyset$ and this concludes the proof of Theorem \ref{bigthpossiml2}.
\epf

\section{Applications and more typical properties of positive contractions}
The purpose of this subsection is to apply the results of the previous section to the study of typical properties of positive contractions on $\ell_p$-spaces for $p >1$, and for the topologies \sot\,, \sott\, and \wot. 

\smallskip

The fact that any Polish topology on $\ppl$ lying between the topologies \wot\, and \sote\, is similar to \sote\, is a very interesting result for the following reason. If the set $\{ T \in \mathcal{P}_1(\ell_2) : T \; \textrm{satisfies the property}\; (P) \}$ is comeager in $(\mathcal{P}_1(\ell_2), \sote)$, then the set $\{ T \in \mathcal{P}_1(\ell_2) : T^* \; \textrm{satisfies the property}\; (P) \}$ is also comeager in $(\mathcal{P}_1(\ell_2), \sote)$. This follows from the fact that the map $T \mapsto T^*$ is a homeomorphism from $(\mathcal{P}_1(\ell_2), \sote)$ onto $(\mathcal{P}_1(\ell_2), \sote)$. This fact remains true on $\ppl$ if we replace the topology \sote\, by the topology \sot, thanks to Theorem \ref{bigthpossiml2}. More generally, we have the following useful fact.

\smallskip

\begin{fact} \label{proptypdualise}

Let $p > 1$ and let (P) be a property of operators on $\ell_p$-spaces.

If a typical $T \in (\mathcal{P}_1(\ell_p), \sot)$ satisfies (P), then a typical $T \in (\mathcal{P}_1(\ell_q),\sott)$ is such that $T^*$ satisfies (P), where $q$ is the conjugate exponent of $p$. This remains true if we replace the topologies \sot\, and \sott\, both by the topology \sote.
\end{fact}

\smallskip
\bpf
This follows from the fact that the map $T \mapsto T^*$ is a homeomorphism from $(\mathcal{P}_1(\ell_p), \sot)$ onto $(\mathcal{P}_1(\ell_q), \sott)$.
\epf

We used this argument in \cite[Proposition 3.2]{Gill} to assert that a typical positive contraction on $(\mathcal{P}_1(\ell_p), \sote)$ is such that $2T^*$ is hypercyclic, since a typical positive contraction $(\mathcal{P}_1(\ell_p), \sote)$ is such that $2T$ is hypercyclic. In particular this was useful to find the point spectrum of a typical positive contraction on $(\mathcal{P}_1(\ell_p), \sote)$, which is empty (we refer to \cite{BM} and \cite{GEP} for background on hypercyclicity).

Also we mention the fact that studying typical properties of positive contractions (resp. of contractions) on $\mathcal{P}_1(\ell_p)$ (resp. on $\mathcal{B}_1(\ell_p)$) for the \wot\, topology is something not always easy. It is also one of the reasons that the notion of similar topologies was studied in \cite{GMM2}. Thanks to Theorems \ref{mainthmposp>1} and \ref{bigthpossiml2}, we will be able to deduce some typical properties of positive contractions on $\ell_p$-spaces for the \wot\,, the \sot\, and the \sott\, topologies.

\subsection{Applications of Theorem \ref{bigthpossiml2}}

The following consequence follows from Theorem \ref{bigthpossiml2}. It concerns the case $p = 2$ and link the comeager sets of all the Polish operator topologies \wot\,, \sott\,, \sot\, and \sote\, on $\ppl$.

\smallskip

\bco \label{corcomeagposp=2}
Let $p = 2$ and let $\tau$ be a Polish operator topology on $\ppl$ lying between the topologies \emph{\wot}\, and \emph{\sote}. Then the topologies $\tau$ and \emph{\sote}\, have the same comeager sets in $\ppl$.
\eco

\smallskip

Corollary \ref{corcomeagposp=2} is far from being true in $\mathcal{B}_1(\ell_2)$ as we explained in subsection \ref{mainresults}. 

The main goal of this article was to give an answer to the open question \cite[Question 5.2]{Gill}, that is to determine the point spectrum of a typical positive contraction $T \in (\ppl, \sot)$. As noted earlier, the point spectrum of a typical positive contraction $T \in (\ppl, \sote)$ is empty because a typical positive contraction $T \in (\ppl, \sote)$ is such that $2T $ and $2T^*$ are hypercyclic. Moreover, there are not many results on the existence of an eigenvalue for positive operators on $\ell_p$-spaces. One of the few relevant results on this subject is \cite[Theorem 5.7]{RT}, which states that any positive compact operator on $\ell_p$ with no non-trivial closed invariant ideals has a unique (up to scaling) positive eigenvector. Although a typical positive contraction $T \in (\ppl, \sot)$ has no non-trivial closed invariant ideals (from \cite[Proposition 4.7]{Gill}), a typical positive contraction $T \in (\ppl, \sot)$ is not compact since a typical positive contraction $T \in (\ppl, \sot)$ is such that $2T$ is hypercyclic. Thanks to Corollary \ref{corcomeagposp=2}, we obtain the following result and we give a negative answer to \cite[Question 5.2]{Gill}. 

\smallskip

\bth \label{corvpspostopotau}
Let $\tau$ be any Polish operator topology on $\ppl$ lying between \emph{\wot}\, and \emph{\sote}. A typical positive contraction $T \in (\ppl, \tau)$ is such that $T$ and $T^*$ have no eigenvalue.
\eth

\smallskip

\bpf
By \cite[Corollary 3.3]{Gill}, a typical positive contraction $T \in (\ppl, \sote)$ is such that $T$ and $T^*$ have no eigenvalue. By Corollary \ref{corcomeagposp=2}, the conclusion of Theorem \ref{corvpspostopotau} immediately follows.
\epf

\smallskip

Another consequence of Corollary \ref{corcomeagposp=2} is the existence of a non-trivial invariant subspace for a typical positive contraction on $\ell_2$.

\smallskip

\bth \label{typinvsubspaceposl2}
Let $\tau$ be any Polish operator topology on $\ppl$ lying between \emph{\wot}\, and \emph{\sote}. A typical positive contraction $T \in (\ppl, \tau)$ has a non-trivial invariant subspace.
\eth

\smallskip

\bpf
By \cite[Proposition 3.1]{Gill}, a typical positive contraction $T \in (\ppl, \sote)$ is such that $\sigma(T) = \overline{\D}$. Thus by Corollary \ref{corcomeagposp=2}, a typical positive contraction $T \in (\ppl, \tau)$ is also such that $\sigma(T) = \overline{\D}$. Now an important result from Brown, Chevreau and Pearcy states that any contraction on a Hilbert space whose spectrum contains the unit circle has a non-trivial invariant subspace (see for instance \cite{BCP2}). Thus the conclusion of Theorem \ref{typinvsubspaceposl2} immediately follows. 
\epf

\subsection{Applications of Theorem \ref{mainthmposp>1}}

The following consequence concerns the cases $1 < p < 2 $ and $p > 2$ and follows from Theorem \ref{mainthmposp>1}.

\smallskip

\bco \label{corcomeagpossotsottsotewot}
Let $X = \ell_p$ with $p > 1$. If $1 < p < 2$, the topologies $\emph{\wot}$ and $\emph{\sot}$ have the same comeager sets in $\ppx$ and the topologies $\emph{\sott}$ and $\emph{\sote}$ have the same comeager sets in $\ppx$. If $p > 2$, the topologies $\emph{\sot}$ and $\emph{\sote}$ have the same comeager sets in $\ppx$ and the topologies $\emph{\wot}$ and $\emph{\sott}$ have the same comeager sets in $\ppx$.
\eco

\smallskip

Using Corollary \ref{corcomeagpossotsottsotewot}, we can describe the point spectrum of a typical positive contraction on $\ell_p$ with $p > 2$ for the topologies \wot\, and\, \sott. This description is analogous to the description given in \cite[Corollary 3.6]{GMM2}.

\smallskip

\bpr \label{comeagvalpropposp>2}
Let $p > 2$ and let $X = \ell_p$. A typical $T \in (\ppx, \emph{\wot})$ (resp. $T \in (\ppx, \emph{\sott})$) has no eigenvalue.
\epr

\smallskip

\bpf
By \cite[Proposition 3.2]{Gill}, we know that a typical $T \in (\mathcal{P}_1(\ell_p), \sot)$ is such that $2T$ is hypercyclic when $1 < p < 2$ and so a typical $T \in (\mathcal{P}_1(\ell_p), \sot)$ is such that $T^*$ has no eigenvalue when $1 < p < 2$. By Fact \ref{proptypdualise}, it follows that a typical $T \in (\mathcal{P}_1(\ell_p), \sott)$ is such that $T$ has no eigenvalue when $p > 2$. Since the topologies \wot\, and \sott\, have the same comeager sets in $\mathcal{P}_1(\ell_p)$ when $p > 2$, Proposition \ref{comeagvalpropposp>2} follows.
\epf

\smallskip

The following result concerns the essential spectrum of a typical positive contraction and is also analogous to \cite[Corollary 3.7]{GMM2}.

\smallskip

\bpr \label{proptypspecess}
Let $p > 1$ and let $X = \ell_p$. A typical $T \in (\ppx, \emph{\wot})$ (resp. $T \in (\ppx, \emph{\sott})$) is such that $\sigma_{ess}(T) = \overline{\D}$, and so $\sigma(T) = \overline{\D}$.
\epr

\smallskip

\bpf
By \cite[Proposition 3.10]{Gill}, we know that a typical $T \in (\ppx, \sot)$ (resp. $T \in (\ppx, \sote)$) is such that $\sigma_{ess}(T) = \overline{\D}$ when $p > 1$. But the property "$\sigma_{ess}(T) = \overline{\D}$" is self-adjoint, so using the same arguments as in the proof of Proposition \ref{comeagvalpropposp>2} and using Corollaries \ref{corcomeagpossotsottsotewot} and \ref{corcomeagposp=2}, we obtain Proposition \ref{proptypspecess}.
\epf

\subsection{Other typical properties of positive contractions for the topologies \wot\, and \sott}

We now focus on other \wot- and \sott-typical properties that have been studied in \cite{Gill} for the \sot\, and the \sote\, topologies. The first property follows from \cite[Proposition 3.7]{Gill}.

\smallskip

\bpr \label{proptypccoisowotsott}
Let $X = \ell_p$ with $p > 1$. A typical $T \in (\ppx, \emph{\wot})$ (resp. $T \in (\ppx, \emph{\sott})$) is such that $T^*$ is not an isometry.
\epr

\smallskip

\bpf
As in the proof of \cite[Proposition 3.7]{Gill}, the set 
$$\mathcal{A} := \displaystyle \bigcup_{j \geq 0} \{ T \in \ppx : \lan e_0^*, T e_j \ran \, \lan e_1^*, T e_j \ran > 0\}$$ 
is contained in the set of positive contractions that are not co-isometries of $X$. The set $\mathcal{A}$ is easily seen to be a \wot-$G_\delta$ and hence a \sott-$G_\delta$ of $\ppx$. Moreover, the set $\mathcal{A}$ is \sote-dense in $\ppx$ and so it is \wot-dense and \sott-dense in $\ppx$. This proves Proposition \ref{proptypccoisowotsott}.
\epf

\smallskip

Recall that an ideal of $X = \ell_p$ is a vector subspace $V$ of $X$ such that if $\lvert x \rvert \leq \lvert y \rvert$ and $y \in V$ then $x \in V$, for every $x,y \in X$. One can show that the closed ideals of $X = \ell_p$ are exactly the sets $[e_k : k \in A]$, where $A$ is a subset of $\Z_+$. As mentioned in \cite{Gill}, there exists a criterion to determine whether a positive operator on $X = \ell_p$ has a non-trivial closed invariant ideal or not. More precisely by \cite[Proposition 1.2]{RT}, a positive operator $T$ on $X = \ell_p$ has no non-trivial closed invariant ideals if and only if
$$
\forall \, i \ne j \geq 0, \, \exists \, n \geq 0 \quad \textrm{such that} \quad \lan e_j^*, T^n e_i \ran > 0.
$$
From this criterion, we deduce the following result.

\smallskip

\bpr \label{proptypidealinvwotsott}
Let $X = \ell_p$ with $p > 1$. A typical $T \in (\ppx, \emph{\wot})$ (resp. $T \in (\ppx, \emph{\sott})$) has no non-trivial closed invariant ideals.
\epr

\smallskip

\bpf
Let us denote by $\mathcal{G}$ the set of positive contractions on $X$ having no non-trivial closed invariant ideals. The proof follows from the fact that the set
$$
\displaystyle \bigcap_{\substack{i, j \geq 0 \\ i \ne j}} \{ T \in \ppx : \lan e_j^*, T e_i \ran > 0 \}
$$
is \wot-$G_\delta$ and hence \sott-$G_\delta$ in $\ppx$, is contained in $\mathcal{G}$ and is \sote-dense in $\ppx$ (see \cite[Proposition 4.7]{Gill} for more details). 
\epf

\smallskip

As mentioned in \cite{Gill}, the existence of non-trivial invariant subspaces for positive operators is a delicate matter. Moreover, it is still unknown whether there exists a positive operator on $\ell_p$ without non-trivial invariant subspace for some $p \geq 1$. The most popular result regarding the existence of non-trivial invariant subspaces for positive operators is the following theorem, due to Abramovich, Aliprantis and Burkinshaw.

\smallskip

\bth[{\cite[Theorem 2.2]{AAB}} {\normalfont and} {\cite[Theorem 2.2]{AAB2}}] \label{thAAB}
Let $X$ be a Banach space with a basis and let $T$ be a positive operator on $X$. Suppose that there exists a non-zero positive operator on $X$ in the commutant of $T$ which is quasinilpotent at a certain non-zero positive vector of $X$. Then $T$ has a non-trivial invariant subspace. Moreover if $X = \ell_p$ and if $p \geq 1$, then $T$ has a non-trivial closed invariant ideal.
\eth

\smallskip

As in \cite{Gill}, we say that a positive operator $T$ on $X = \ell_p$ satisfies the AAB criterion if it satisfies the hypotheses of Theorem \ref{thAAB}. From Proposition \ref{proptypidealinvwotsott}, the following result immediately follows.

\smallskip

\bpr \label{proptypAAB}
Let $X = \ell_p$ with $p > 1$. A typical $T \in (\ppx, \emph{\wot})$ (resp. $T \in (\ppx, \emph{\sott})$) does not satisfy the AAB criterion.
\epr


\section{More on the points of continuity of the identity map on $\mathcal{P}_1(\ell_p)$}\label{Section3}
The aim of this section is to study in more depth the sets $\mathcal{C}(\wot, \sot)$, $\mathcal{C}(\wot, \sott)$ and $\mathcal{C}(\sot, \sote)$. We will give examples of such points of continuity when $X = \ell_2$. We will also show that, when $X = \ell_2$, these sets are very different from the sets of the points of continuity of the identity maps on $\mathcal{B}_1(\ell_2)$ for the corresponding topologies. 

\smallskip

As mentioned in Subsection \ref{subsect2a}, the set $\mathcal{C}(\sot, \sote)$ is comeager in $(\mathcal{P}_1(\ell_p), \sot)$ when $p > 1$, while the sets $\mathcal{C}(\wot, \sot)$ and $\mathcal{C}(\wot, \sott)$ are comeager in $(\mathcal{P}_1(\ell_p), \wot)$ when $p > 1$. So the description of these sets can be used to determine new typical properties of positive contractions on $\ell_p$-spaces. 

\smallskip

In \cite[Proposition 2.11]{GMM2}, it was proved that the set of all points of continuity of the identity map from $(\mathcal{B}_1(\ell_2), \sot)$ into $(\mathcal{B}_1(\ell_2), \sote)$ (resp. from $(\mathcal{B}_1(\ell_2), \wot)$ into $(\mathcal{B}_1(\ell_2), \sott)$) is the co-isometries of $\ell_2$. It was also proved that the set of all points of continuity of the identity map from $(\mathcal{B}_1(\ell_2), \wot)$ into $(\mathcal{B}_1(\ell_2), \sot)$ is the isometries of $\ell_2$. These results are no longer true for positive contractions. Indeed by \cite[Proposition 3.7]{Gill}, a typical positive contraction $T \in (\ppl, \sot)$ is not a co-isometry, and also a typical positive contraction $T \in (\ppl, \wot)$ is not an isometry for the same reasons. We will see in this section that the situation is more complicated in the positive setting.

\smallskip

We first provide some observations concerning the sets $\mathcal{C}(\wot, \sott)$, $\mathcal{C}(\sot, \sote)$ and $\mathcal{C}(\wot, \sot)$. If $p > 1$, we identify $(\ell_p)^*$ with $\ell_q$, where $q$ is the conjugate exponent of $p$. 

\smallskip

\bpr \label{proprmkptscont}
Let $X = \ell_p$ with $p > 1$ and let $T$ be a positive contraction on $X$. Then $T \in \mathcal{C}(\emph{\wot}, \emph{\sott})$ if and only if $T^* \in \mathcal{C}(\emph{\wot}, \emph{\sot})$.
\epr

\smallskip

\bpf
Let $T$ be a positive contraction on $X$. Suppose that $T$ belongs to $\mathcal{C}(\wot, \sott)$ and let us show that $T^* \in \mathcal{C}(\wot, \sot)$. Let $(T_n)_{n \geq 0}$ be a sequence of positive contractions such that $T_n \underset{n \to \infty}{\overset{\wot}{\longrightarrow}}T^*$. Then $T_n^* \underset{n \to \infty}{\overset{\wot}{\longrightarrow}} T$.
Thus $T_n^* \underset{n \to \infty}{\overset{\sott}{\longrightarrow}} T$ since $T \in \mathcal{C}(\wot, \sott)$, and so $T_n \underset{n \to \infty}{\overset{\sot}{\longrightarrow}}T^*$. This proves that $T^* \in \mathcal{C}(\wot, \sot)$. The same arguments show that the converse is true and this concludes the proof of Proposition \ref{proprmkptscont}. 
\epf

We now discuss some elementary properties that an element of $\mathcal{C}(\sot, \sote)$, $\mathcal{C}(\wot, \sott)$ or $\mathcal{C}(\wot, \sot)$ must have. The first ones concern the norm of a point of continuity.

\smallskip

\bpr \label{ptcontnormsotsote}
Let $T \in \mathcal{P}_1(\ell_2)$. If $T$ belongs to one of the sets\, $ \mathcal{C}(\emph{\sot}, \emph{\sote}),\, \mathcal{C}(\emph{\wot}, \emph{\sot})$ or $\mathcal{C}(\emph{\wot}, \emph{\sott})$ , then $\lVert T \rVert = 1$.
\epr

\smallskip

\bpf
By Proposition \ref{proprmkptscont}, it is enough to prove that if $T$ belongs to $\mathcal{C}(\sot, \sote)$, then $\lVert T \rVert = 1$. Suppose that $\lVert T \rVert < 1$ and let $\delta > 0$ be such that $\lVert T \rVert + \delta < 1$.
Let us consider the positive operator $T_n$ on $\ell_2$ defined by $T_n x = T P_n x + \delta \, \lan e_{n+1}, x \ran \, e_0$ for every $x \in \ell_2$ and every $n \geq 0$. Our choice of $\delta
$ implies that each $T_n$ is a positive contraction on $\ell_2$. We also have $T_n \underset{n \to \infty}{\overset{\sot}{\longrightarrow}} T$. Moreover 
$$
\lVert T_n ^* e_0 - T^* e_0 \rVert^2 \geq (\delta - \lan e_{n+1}, T^* e_0 \ran)^2,
$$
so $\displaystyle \liminf_{n \to \infty} \lVert T_n ^* e_0 - T^* e_0 \rVert^2 \geq \delta^2 > 0 $ and $T_n \notin \mathcal{C}(\sot, \sote)$. This proves Proposition \ref{ptcontnormsotsote}.
\epf

\smallskip
Our second set of observations concern properties of the rows and the lines of the matrices representing points of continuity (in the canonical basis). 

\smallskip

\bpr \label{ptcontlignes}
Let $T \in \mathcal{P}_1(\ell_2)$. If $T \in \mathcal{C}(\emph{\sot}, \emph{\sote})$, then $T^* e_l \ne 0$ for every $l \geq 0$.
\epr

\smallskip

\bpf
Suppose that there exists an index $l \geq 0$ such that $T^* e_l = 0$. Then $T = P_{l-1} T + Q_l T$. Let us consider the positive operator $T_n$ on $\ell_2$ defined by $T_n x = T P_n x + \lan x , e_{n+1} \ran \, e_l$ for every $x \in \ell_2$ and every $n \geq 0$. Each operator $T_n$ is a contraction since for every $x \in \ell_2$, we have 
\begin{align*}
    \lVert T_n^* x \rVert^2 &= \lVert P_n T^* P_{l-1} x + \lan x, e_l \ran\, e_{n+1} + P_n T^* Q_l x \rVert^2 \\
    &= \lvert \lan x, e_l \ran \rvert^2 + \lVert P_n T^* (P_{l-1} x + Q_l x) \rVert^2 \\
    &\leq \lvert \lan x, e_l \ran \rvert^2 + \lVert P_n T^* \rVert^2 \lVert P_{l-1} x + Q_l x \rVert^2 \\
    &\leq \lvert \lan x, e_l \ran \rvert^2 + \displaystyle \sum_{\substack{j \geq 0 \\j \ne l}} \lvert \lan x, e_j \ran \rvert^2 \\
    &\leq \lVert x \rVert^2.
\end{align*}
As in the proof of Proposition \ref{ptcontnormsotsote}, the sequence $(T_n)$ converges to $T$ for the \sot\, topology but does not converge to $T$ for the \sote\, topology. This implies that $T \notin \mathcal{C}(\sot, \sote)$.
\epf

\smallskip

From Propositions \ref{proprmkptscont} and \ref{ptcontlignes}, we deduce the following one.

\smallskip

\bpr \label{ptcontcol}
Let $T \in \mathcal{P}_1(\ell_2)$. If $T \in \mathcal{C}(\emph{\wot}, \emph{\sott}$) (resp. $T \in \mathcal{C}(\emph{\wot}, \emph{\sot})$), then $T^* e_l \ne 0$ for every $l \geq 0$ (resp. $T e_l \ne 0$ for every $l \geq 0$).
\epr

\smallskip

Any co-isometry of $\ell_2$ belongs to $\mathcal{C}(\wot, \sott)$ and any isometry of $\ell_2$ belongs to $\mathcal{C}(\wot, \sot)$ (the same arguments as in \cite[Proposition 2.11]{GMM2} hold). 

\smallskip

Our aim is now to exhibit new sets of points of continuity. To this aim, we need the following lemma.

\smallskip

\blm \label{lemptscontwotsott}
Let $B \in \mathcal{P}_1(\ell_2)$. Then $B$ belongs to $\mathcal{C}(\emph{\wot}, \emph{\sott})$ if and only if the following property holds: for every $\varepsilon > 0$ and for every $r \geq 0$, there exist $\delta > 0$ and $m \geq 0$ such that  
\begin{align} \label{carsotsotecont}
    \textrm{if} \; T \in \mathcal{P}_1(\ell_2) \; &\textrm{is such that} \; \displaystyle \max_{0 \leq k,l \leq m} \lvert \lan e_l, (T-B) e_k \ran \rvert < \delta, \\ \notag
    &\textrm{then} \; \displaystyle \max_{0 \leq k \leq r} \lVert (T-B)^* e_k \rVert < \varepsilon.
\end{align}
\elm

\smallskip
The proof of Lemma \ref{lemptscontwotsott} uses similar arguments to the proof of \cite[Lemma 2.2]{Gill}, that is uses local bases for the topologies \wot\, and \sott\, and the triangular inequality. For this reason, we will omit its proof. We now introduce two new classes of points of continuity, denoted by $\mathcal{M}$ and $\mathcal{M}'$. The set $\mathcal{M}$ strictly contains the co-isometries of $\ell_2$, whereas the set $\mathcal{M}'$ strictly contains the isometries of $\ell_2$. 

\smallskip

\bth \label{thptscontsotsoteunionsuppN}
Let $\mathcal{M}$ be the class of every positive contractions $T$ on $\ell_2$ such that
\begin{itemize}
    \item $\lVert T \rVert = 1$;
    \item there exists a family $(u_r)_{r \in J}$ of norming vectors for $T^*$ indexed by a set $J \subseteq \Z_+$ such that $u_r \geq 0$, $\lVert u_r \rVert = 1$ for every $r \in J$ and $\displaystyle \bigcup_{r \in J} \textrm{Supp}(u_r) = \Z_+ $.
\end{itemize}
Then every operator in $\mathcal{M}$ belongs to $\mathcal{C}(\emph{\wot}, \emph{\sott})$ and to $\mathcal{C}(\emph{\sot}, \emph{\sote})$.
\eth

\smallskip

\bpf
The proof of Theorem \ref{thptscontsotsoteunionsuppN} is not very far from the proof of Proposition \ref{lem3.3pos2}. The only difference is that here, the operators of the class $\mathcal{M}$ and the vectors $u_r$ are defined on $\ell_2$ and not just on $E_M$ for some $M \geq 1$. We must therefore pay attention to the fact that some vectors $u_r$ can have an infinite support. 

\smallskip

Let $B \in \mathcal{M}$.
For every $k \in J$ and for every $l \geq 0$, we set for this proof $u_k(l) := \lan e_l, u_k \ran$. Also recall that $b_{i,j} = \lan e_i, Be_j \ran$ for $i,j \geq 0$. Let $\varepsilon > 0$ and let $r \geq 0$. We have to find $\delta > 0$ and an integer $n \geq 0$ such that (\ref{carsotsotecont}) holds.

\smallskip

By the condition on the supports, we can write $\{0,...,r \} = \displaystyle \bigcup_{k \in J_r} (\textrm{Supp}(u_k) \cap \{ 0,...,r \})$, where $J_r$ is a certain finite subset of $J$. Let us now consider the following two sets
\begin{align*}
    I := \{ k \in J_r : \; \textrm{Supp}(u_k) \; \textrm{is infinite}\} \quad \textrm{and} \quad I' := \{ k \in J_r : \; \textrm{Supp}(u_k) \; \textrm{is finite} \}.
\end{align*}
For $k \in I'$, let us write $n_k = \max \textrm{Supp}(u_k)$. For every $k \in I'$, one can find an integer $n_k' > r$ such that
\begin{align}
   \label{eq1 thsotsoteunionsupp} \displaystyle &\max_{\substack{0 \leq l \leq r \\ u_k(l) > 0}} \frac{1}{(u_k(l))^2} \, (1- \lVert P_{n_k'} B^* u_k  \rVert^2) < \varepsilon/4 \\
   \label{eq2 thsotsoteunionsupp} \textrm{and} \quad &\displaystyle \max_{0 \leq l \leq r} \lVert Q_{n_k'} B^* e_l \rVert^2 < \varepsilon/4.
\end{align}
Let us write $n =\displaystyle \max_{k \in I'} n_k$ and $n' = \displaystyle \max_{k \in I'} n_k'$.

\smallskip

Similarly for every $k \in I$, one can find an integer $N_k > r$ such that
\begin{align}
   \label{e3 thsotsoteunionsupp} \displaystyle &\max_{\substack{0 \leq l \leq r \\ u_k(l) > 0}} \frac{1}{(u_k(l))^2} \,( 1 -  \lVert P_{N_k} B^* P_{N_k} u_k  \rVert^2) < \varepsilon/4 \\
   \label{eq4 thsotsoteunionsupp} \textrm{and} \quad &\displaystyle \max_{0 \leq l \leq r} \lVert Q_{N_k} B^* e_l \rVert^2 < \varepsilon/4.
\end{align}
We also write $N = \displaystyle \max_{k \in I} N_k$. 

\smallskip

We choose $\delta > 0$ such that $b_{i,j} - \delta > 0$ for every $0 \leq i,j \leq \max(n, N, n') + 1$ such that $b_{i,j} > 0$, 
\begin{align}
    2 \delta &\displaystyle \max_{\substack{0 \leq l \leq r \\ k \in J_r \\ u_k(l) > 0}} \frac{1}{(u_k(l))^2} \left(\displaystyle \sum_{j=0}^n (u_k(j))^2 \sum_{m=0}^{n'} b_{j,m} + \displaystyle \sum_{0 \leq i < j \leq n} u_k(i) u_k(j) \sum_{m=0}^{n'} (b_{i,m} + b_{j,m}) \right)\\ \notag
    &< \varepsilon/2
\end{align}
and
\begin{align}
   2 \delta &\displaystyle \max_{\substack{0 \leq l \leq r \\ k \in J_r \\ u_k(l) > 0}} \frac{1}{(u_k(l))^2} \left(\displaystyle \sum_{j=0}^N (u_k(j))^2 \sum_{m=0}^{N} b_{j,m} + \displaystyle \sum_{0 \leq i < j \leq N} u_k(i) u_k(j) \sum_{m=0}^{N} (b_{i,m} + b_{j,m}) \right)\\\notag
   &< \varepsilon / 2. 
\end{align}

\smallskip

Let $T \in \ppl$ be such that $$\displaystyle \max_{0 \leq k, l \leq \max(n, N, n') + 1} \lvert \lan e_l, (T-B) e_k \ran \rvert < \delta. $$ Let us show that $\displaystyle \max_{0 \leq k \leq r} \lVert (T-B)^* e_k \rVert^2 < \varepsilon$.

\smallskip

We will control $\lVert (T-B)^* e_l \rVert^2$ for every $l \in \textrm{Supp}(u_k) \cap \{ 0 , ... , r \}$, when $k \in I$ and when $k \in I'$.  

\smallskip

For every $k \in I'$, one has
\begin{align}
    \label{th36eq1} 1 \geq \Big \lVert T^* \Big(\displaystyle \sum_{j=0}^{n_k} u_k(j) e_j \Big) \Big \rVert^2 = &\displaystyle \sum_{j=0}^{n_k} (u_k(j))^2 \lVert T^* e_j \rVert^2 \\ \notag
     &+ 2 \displaystyle \sum_{0 \leq i < j \leq n_k} u_k(i) u_k(j) \, \lan T^* e_i, T^* e_j \ran.  
\end{align}
Using positivity and the fact that $$\displaystyle \max_{0 \leq k,l \leq \max(n,N,n') + 1} \lvert \lan e_l, (T-B) e_k \ran \rvert < \delta,$$ we can say that
\begin{align}
    \label{th36eq2} \lan T^* e_i, T^* e_j \ran \geq \displaystyle \sum_{\substack{0 \leq m \leq n_k' \\ b_{i,m} > 0, b_{j,m} > 0}} (b_{i,m} - \delta) (b_{j,m} - \delta)
\end{align}
for every $0 \leq i,j \leq n_k$.

\smallskip

Using (\ref{th36eq1}), (\ref{th36eq2}) and the positivity of the vectors $u_k$, we thus have for every $0 \leq l \leq n_k$ and for every $k \in I'$
\smallskip

\begin{align}
  \label{th36eq3}  (u_k(l))^2 \lVert T^* e_l \rVert^2 &\leq 1 - \displaystyle \sum_{\substack{0 \leq j \leq n_k \\ j \ne l}} (u_k(j))^2 \displaystyle \sum_{\substack{0 \leq m \leq n_k' \\ b_{j,m} > 0}} (b_{j,m} - \delta)^2 \\ \notag 
    &- 2 \displaystyle \sum_{0 \leq i < j \leq n_k} u_k(i) u_k(j) \sum_{\substack{0 \leq m  \leq n_k' \\b_{i,m} > 0, \, b_{j,m} > 0}} (b_{i,m} - \delta) (b_{j,m} - \delta).
\end{align}

\smallskip
As (\ref{th36eq2}), one can say that
\begin{align}
   \label{th36eq4} \lan T^* e_l, B^* e_l \ran \geq \displaystyle \sum_{\substack{0 \leq m \leq n_k' \\ b_{l,m} > 0}} b_{l,m} (b_{l,m} - \delta)
\end{align}
for every $l \in \textrm{Supp}(u_k) \cap \{0,...,r \}$ such that $k \in I'$.

\smallskip

Thus using (\ref{th36eq3}) and (\ref{th36eq4}), the quantity 
\begin{align}
  \label{th36eq5}  &(u_k(l))^2 \lVert T^* e_l - B^* e_l \rVert^2 \\ \notag
  &= (u_k(l))^2 \lVert T^* e_l \rVert^2 - 2 (u_k(l))^2 \lan T^* e_l , B^* e_l \ran + (u_k(l))^2 \lVert B^* e_l \rVert^2
\end{align}
can be bounded above for every $k \in I'$ and for every $l \in \textrm{Supp}(u_k) \cap \{ 0,...,r \}$ by
\begin{align}
   \label{th36eq6} &1 - \displaystyle \sum_{\substack{0 \leq j \leq n_k \\ j \ne l}} (u_k(j))^2 \displaystyle \sum_{\substack{0 \leq m \leq n_k' \\ b_{j,m} > 0}} (b_{j,m} - \delta)^2\\ \notag 
   &- 2 \displaystyle \sum_{0 \leq i < j \leq n_k} u_k(i) u_k(j) \sum_{\substack{0 \leq m  \leq n_k' \\b_{i,m} > 0, \, b_{j,m} > 0}} (b_{i,m} - \delta) (b_{j,m} - \delta) \\ \notag
    & -2 (u_k(l))^2 \displaystyle \sum_{\substack{0 \leq m \leq n_k' \\ b_{l,m} > 0}} b_{l,m} (b_{l,m} - \delta) + (u_k(l))^2 \displaystyle \sum_{m \geq 0} b_{l,m}^2\, , 
\end{align}
that is by
\begin{align}
  \label{th36eq7}  1 &- \left( \displaystyle \sum_{j=0}^{n_k} (u_k(j))^2 \displaystyle \sum_{m=0}^{n_k'} b_{j,m}^2 + 2 \displaystyle \sum_{0 \leq i < j \leq n_k} u_k(i) u_k(j) \displaystyle \sum_{m=0}^{n_k'} b_{i,m} b_{j,m} \right) \\ \notag
    + &2 \delta \left(\displaystyle \sum_{j=0}^n (u_k(j))^2 \sum_{m=0}^{n'} b_{j,m} + \displaystyle \sum_{0 \leq i < j \leq n} u_k(i) u_k(j) \sum_{m=0}^{n'} (b_{i,m} + b_{j,m}) \right)\\ \notag 
    &+ (u_k(l))^2 \lVert Q_{n_k'} B^* e_l \rVert^2
\end{align}
since $n \geq n_k$ and $n' \geq n_k'$.

\smallskip

Therefore the quantity (\ref{th36eq5}) can be bounded above by
\begin{align}
    \label{th36eq8} &(1 - \lVert P_{n_k'} B^* u_k \rVert^2) \\ \notag
     &+ 2 \delta \left(\displaystyle \sum_{j=0}^n (u_k(j))^2 \sum_{m=0}^{n'} b_{j,m} + \displaystyle \sum_{0 \leq i < j \leq n} u_k(i) u_k(j) \sum_{m=0}^{n'} (b_{i,m} + b_{j,m}) \right) \\ \notag 
    &+ (u_k(l))^2 \lVert Q_{n_k'} B^* e_l \rVert^2,
\end{align}
since 
\begin{align*}
    \displaystyle \sum_{j=0}^{n_k} (u_k(j))^2 \displaystyle \sum_{m=0}^{n_k'} b_{j,m}^2 + 2 \displaystyle \sum_{0 \leq i < j \leq n_k} u_k(i) u_k(j) \displaystyle \sum_{m=0}^{n_k'} b_{i,m} b_{j,m} &= \displaystyle \sum_{m=0}^{n_k'} \lan u_k, P_{n_k} B e_m \ran^2 \\ &= \lVert P_{n_k'} B^* u_k \rVert^2.
\end{align*}

\smallskip

By the same way, replacing $n_k$ and $n_k'$ by $N_k$, we have for every $k \in I$ and for every $l \in \textrm{Supp}(u_k) \cap \{ 0,...,r \}$,
\begin{align}
   \label{th36eq9} & (u_k(l))^2 \lVert T^* e_l - B^* e_l \rVert^2  \\ \notag
   &\leq (1 - \lVert P_{N_k} B^* P_{N_k} u_k \rVert^2) \\ \notag
   &+ 2 \delta \left(\displaystyle \sum_{j=0}^N (u_k(j))^2 \sum_{m=0}^{N} b_{j,m} + \displaystyle \sum_{0 \leq i < j \leq N} u_k(i) u_k(j) \sum_{m=0}^{N} (b_{i,m} + b_{j,m}) \right) \\ \notag 
    &+ (u_k(l))^2 \lVert Q_{N_k} B^* e_l \rVert^2
\end{align}
since $N \geq N_k$.

\smallskip
Using (\ref{eq1 thsotsoteunionsupp}), (\ref{eq2 thsotsoteunionsupp}), (\ref{e3 thsotsoteunionsupp}), (\ref{eq4 thsotsoteunionsupp}), (\ref{th36eq8}) and (\ref{th36eq9}), we easily see that $$\displaystyle \max_{0 \leq l \leq r} \lVert (B-T)^* e_l \rVert^2 < \varepsilon.$$ This concludes the proof of Theorem \ref{thptscontsotsoteunionsuppN}.
\epf

\smallskip

From Proposition \ref{proprmkptscont} and Theorem \ref{thptscontsotsoteunionsuppN}, we obtain the following result.

\smallskip

\bth \label{thptscontwotsotunionsuppN}
Let $\mathcal{M}'$ be the class of every positive contractions $T$ on $\ell_2$ such that
\begin{itemize}
    \item $\lVert T \rVert = 1$;
    \item there exists a family $(u_r)_{r \in J}$ of norming vectors for $T$ indexed by a set $J \subseteq \Z_+$ such that $u_r \geq 0$, $\lVert u_r \rVert = 1$ for every $r \in J$ and $\displaystyle \bigcup_{r \in J} \textrm{Supp}(u_r) = \Z_+ $.
\end{itemize}
Then every operator in $\mathcal{M}'$ belongs to $\mathcal{C}(\emph{\wot}, \emph{\sot})$.
\eth

\smallskip

We will now see that the classes $\mathcal{M}$ and $\mathcal{M}'$ are \sot-dense in $\ppl$. They thus are large classes of continuity points.
To do so, we will use the following approximation lemma, which is proved in \cite[Lemma 2.2]{Gill}.

\smallskip

\blm \label{lemapprox}
Let $X = \ell_p$ with $p > 1$. Let $\mathcal{C}(X)$ be a class of operators on $X$ and define $\mathcal{C}_1(X) := \mathcal{C}(X) \cap \ppx$.

Suppose that the following property holds: for every $\varepsilon >0,$ every $N \in \Z_+$ and for every positive operator $A \in \mathcal{P}(E_N)$ with $\lVert A \rVert < 1$, there exists a positive operator $T \in \mathcal{C}_1(X)$ such that
\begin{equation} \label{eq1lemapprox}
    \lVert (T-A) \, e_k \rVert < \varepsilon \quad \textrm{for every} \quad  0 \leq k \leq N.
\end{equation}
Then $\mathcal{C}_1(X)$ is dense in $(\ppx, \emph{\sot})$.
\elm

\smallskip

With Lemma \ref{lemapprox}, we obtain the following result.

\smallskip

\bpr \label{proplargeclasses}
The classes $\mathcal{M}$ and $\mathcal{M}'$ are \emph{\sot}-dense in $\ppl$.
\epr

\smallskip

\bpf
In order to apply Lemma \ref{lemapprox}, let $N \geq 0$, let $\varepsilon > 0$ and let $A \in \mathcal{P}_1(E_N)$ be such that $\lVert A \rVert < 1$. By Proposition \ref{introprop5.15}, one can find a positive operator $B$ on $E_{2N+1}$ such that $\lVert B \rVert = 1$, $B$ has a positive norming vector $u \in E_{2N+1}$ for which $\textrm{Supp}(u) = \{ 0,..., 2N+1 \}$ and $\lVert B P_N - A \rVert < \varepsilon$. Let us define a positive operator $T$ on $\ell_2$ by $T = B P_{2N+1} + Q_{2N +1}$. Then $\lVert T \rVert = 1$ and, since the operator $Q_{2N+1}$ has a positive norming vector on $F_{2N+1}$ with positive coordinates, we can easily obtain a positive norming vector for $T$ whose support is $\Z_+$. This shows that $T$ belongs to $\mathcal{M}'$. Moreover, we have $\displaystyle \max_{0 \leq k \leq N} \lVert (T-A) e_k \rVert < \varepsilon$. This shows that $\mathcal{M}'$ is \sot-dense in $\ppl$. Also notice that by Proposition \ref{introprop5.15} and by Remark \ref{remarkvectnormant}, the operator $B^*$ also has a positive norming vector whose support is $\{ 0,..., 2N+1 \}$ and so by the same way, the operator $T $ belongs to $\mathcal{M}$. This concludes the proof of Proposition \ref{proplargeclasses}.     
\epf

\smallskip

Finally we end this article by showing that the class $\mathcal{M}$ is strictly contained in $\mathcal{C}(\wot, \sott)$ and in $\mathcal{C}(\sot, \sote)$, and that $\mathcal{M}'$ is strictly contained in $\mathcal{C}(\wot, \sot)$. To do so, we need the following important proposition.

\smallskip

\bpr \label{proptypnormattain}
Let $X = \ell_p$ with $p > 1$. A typical $T \in (\ppx, \emph{\sote})$ does not attain its norm.
\epr

\smallskip

\bpf
Let $\mathcal{N}$ be the set of positive contractions on $X$ which attain their norm. Since the set $S^+(X) := \{ T \in \px : \lVert T \rVert = 1 \}$ is comeager in $(\ppx, \sote)$ by \cite[Proposition 3.1]{Gill}, it is enough to prove that $S^+(X) \setminus\, \mathcal{N}$ is comeager in $(S^+(X), \sote)$. As in the proof of \cite[Theorem 6.1]{GMM2}, the set $S^+(X) \setminus\, \mathcal{N}$ is a $G_\delta$ of $(S^+(X),\sote)$. This set is also \sote-dense in $S^+(X)$. Indeed, let $A \in S^+(X)$ and let $(\varepsilon_n)_{n \geq 0}$ be a sequence of positive real numbers tending to $0$. For every $N \geq 0$, let $B_N$ be a positive operator of norm $1$ on $F_N$ which does not attain its norm on $F_N$ (for example we can take a diagonal operator on $F_N$ associated to a sequence of positive real numbers in $(0,1)$ which converges to $1$). If we set $T_N = (1- \varepsilon_N) P_N A P_N + Q_N B_N Q_N$ for every $N \geq 0$, then each $T_N$ is a positive operator in $S^+(X)$, and we can show as in \cite[Theorem 6.1]{GMM2} that these operators do not attain their norm. Since the sequence $(T_N)$ converges to $A$ for the \sote\, topology, the set $S^+(X) \setminus\, \mathcal{N}$ is \sote-dense in $(S^+(X), \sote)$. This proves Proposition \ref{proptypnormattain}.
\epf

\smallskip

Corollary \ref{corcomeagposp=2} has the following important consequence.

\smallskip

\bco \label{corproptypnormattain}
A typical $T \in (\ppl, \emph{\sot})$ (resp. $T \in (\ppl, \emph{\wot})$, $T \in  (\ppl, \emph{\sott})$) does not attain its norm.
\eco

\smallskip

With Corollary \ref{corcomeagpossotsottsotewot}, we also obtain:

\smallskip

\bco \label{corproptypnormattainlp}
Let $X = \ell_p$ with $p>1$. If $p >2$, a typical $T \in (\ppx, \emph{\sot})$ does not attain its norm. If $1 < p < 2$, a typical $T \in (\ppx, \emph{\sott})$ does not attain its norm.
\eco

\smallskip

Theorem \ref{thptscontsotsoteunionsuppN} (resp. Theorem \ref{thptscontwotsotunionsuppN}) asserts that the class $\mathcal{M}$ (resp. $\mathcal{M}'$) is contained in the sets\, $\mathcal{C}(\wot, \sott)$ and\, $\mathcal{C}(\sot, \sote)$ (resp. in\, $\mathcal{C}(\wot, \sot)$). By Corollary \ref{corproptypnormattain}, we know that a typical positive contraction on $\ell_2$ for the \sot\, and the \wot\, topologies does not attain its norm, so in particular the class $\mathcal{M}$ can't be equal to $\mathcal{C}(\sot, \sote)$ or to $\mathcal{C}(\wot, \sott)$, and the class $\mathcal{M}'$ can't be equal to $\mathcal{C}(\wot, \sot)$. This shows that the description of the sets $\mathcal{C}(\wot, \sott)$, $\mathcal{C}(\sot, \sote)$ and $\mathcal{C}(\wot, \sot)$ in the positive setting may be an arduous problem, contrary to the description of these sets in the general setting of contractions on $\ell_2$.

\section{Comments and open problems}\label{Questions}
We end this article with some comments and open questions related to our previous results.

As mentioned at the end of Section \ref{Section3}, the descriptions of the sets $\mathcal{C}(\wot, \sott)$, $\mathcal{C}(\sot, \sote)$ and $\mathcal{C}(\wot, \sot)$ are not yet complete. So the following problem it still open.

\smallskip

\begin{problem}
    Describe completely the sets $\mathcal{C}(\wot, \sott)$, $\mathcal{C}(\sot, \sote)$ and $\mathcal{C}(\wot, \sot)$.
\end{problem}

\smallskip

So far, we have proved that the topologies \wot, \sott, \sot\, and \sote\, are similar on $\ppl$. Also when $p> 2$, the topologies \sot\, and \sote\, are similar on $\mathcal{P}_1(\ell_p)$ and the topologies \wot\, and \sott\, are similar on $\mathcal{P}_1(\ell_p)$. And finally when $1 < p < 2$, the topologies \sott\, and \sote\, are similar on $\mathcal{P}_1(\ell_p)$ and the topologies \wot\, and \sot\, are similar on $\mathcal{P}_1(\ell_p)$. The proof of Theorem \ref{bigthpossiml2} is very particular to the Hilbertian setting, since it uses the parallelogram identity in the proof of Proposition \ref{lem3.3pos2}. This proof cannot be used for other values of $p$. Thus the following problem is open.

\smallskip

\begin{problem} \label{pb2}
    Is it true that the topologies \wot, \sott, \sot\, and \sote\, are similar on $\mathcal{P}_1(\ell_p)$ for $1 < p < 2$ and for $p > 2$?
\end{problem}

\smallskip

A positive answer to Open Problem \ref{pb2} would not be surprising, and it would also not be surprising if the notion of norming vector helps to solve this problem. In particular if the answer to Open Problem \ref{pb2} is affirmative, then the point spectrum of a typical positive contraction on $\ell_p$ for $p > 1$ would be empty for the topologies \wot, \sot, \sott\, and \sote.

\medskip


\end{document}